\documentclass[11pt]{amsart}
\usepackage{amsfonts,amssymb,eucal}

\newtheorem{thm}{Theorem}[section]
\newtheorem{cor}[thm]{Corollary}

\newtheorem{prop}[thm]{Proposition}

\newcommand{\inprod}[2]{\left\langle #1, #2 \right\rangle}

\DeclareMathOperator{\Var}{Var}

\DeclareMathOperator{\conv}{conv}

\newcommand{\Prob}{\mathbb{P}}
\newcommand{\R}{\mathbb{R}}
\newcommand{\E}{\mathbb{E}}
\renewcommand{\epsilon}{\varepsilon}

\allowdisplaybreaks[3]
\numberwithin{equation}{section}

\title{Gaussian marginals of convex bodies with symmetries}

\author[M.\ Meckes]{Mark W.\ Meckes} 

\email{mark.meckes@case.edu}

\address{Department of Mathematics, Case Western Reserve University,
Cleveland, Ohio 44106, U.S.A.}

\begin{document}

\begin{abstract}
We prove Gaussian approximation theorems for specific $k$-dimensional
marginals of convex bodies which possess certain symmetries. In
particular, we treat bodies which possess a $1$-unconditional basis,
as well as simplices.  Our results extend recent results for
$1$-dimensional marginals due to E.\ Meckes and the author.
\end{abstract}

\maketitle

%%%%%%%%%%%%%%%%%%%%%%%%%%%%%%%%%%%%%%%%
\section{Introduction}\label{S:intro}

Let $K$ be a convex body in the Euclidean space $\R^n$, $n\ge 2$,
equipped with its standard inner product $\inprod{\cdot}{\cdot}$ and
Euclidean norm $|\cdot|$, and let $\mu$ denote the uniform (normalized
Lebesgue) probability measure on $K$. In this paper we consider
$k$-dimensional marginals of $\mu$, that is, the push-forward $\mu
\circ P_E^{-1}$ of $\mu$ by the orthogonal projection $P_E$ onto some
$k$-dimensional subspace $E\subset \R^n$.

The question of whether every convex body $K\subset \R^n$ has
$1$-dimensional marginals which are close to Gaussian measures when
$n$ is large is known as the \emph{central limit problem for convex
bodies}, and was apparently first explicitly posed in the literature
in \cite{ABP,BV}. A natural extension is to ask, for how large $k\le
n$ does $K$ necessarily possess nearly Gaussian $k$-dimensional
marginals? The latter question can be thought of as asking for a
measure-theoretic analogue of Dvoretzky's theorem, which implies the
existence of nearly ellipsoidal $k$-dimensional projections of $K$
when $k \ll \log n$.

Very recently Klartag \cite{Klartag,Klartag2} showed that any convex
body has nearly Gaussian $k$-dimensional marginals when $k\le c
n^\alpha$, where $c>0$ and $0<\alpha<1$ are some universal constants;
closeness of probability measures is quantified by the total variation
metric and also has a power-law dependence on $n$. This points out an
important difference from Dvoretzky's theorem, in which it is known
that for an arbitrary convex body $k$ can only be taken to be
logarithmically large in $n$.  Klartag's work followed partial
results, involving different additional hypotheses and metrics between
probability measures, by many authors; we mention
\cite{BB,Milman,Sodin,Klartag3,EK} among recent contributions and
refer to \cite{Klartag} for further references.

In much of the work on this problem, including the main results of
\cite{Klartag,Klartag2}, the existence of nearly Gaussian marginals
$\mu\circ P_E^{-1}$ is proved nonconstructively, so that no concrete
such subspace $E$ is exhibited. This is typical of the proofs of
Dvoretzky-like results.  In \cite{MM}, E.\ Meckes and the author used
Stein's method of exchangeable pairs to prove Berry-Esseen theorems
for specific $1$-dimensional marginals of convex bodies which possess
certain types of symmetries. Roughly, under some additional
hypotheses, \cite{MM} shows that a $1$-dimensional marginal $\mu\circ
P_E^{-1}$ is nearly Gaussian when $K$ possesses many symmetries
$\sigma:\R^n\to \R^n$ for which the $1$-dimensional subspace
$\sigma(E)\subset \R^n$ is very different from $E$. In
\cite{Klartag,Klartag3}, another approach is used to study the
marginal of a $1$-unconditional body on the subspace spanned by
$(1,\dotsc,1)$; see the remarks at the end of this paper for the
relationship between these approaches.

The main purpose of this paper is to prove versions of the results of
\cite{MM} for $k$-dimensional marginals with $k\ge 2$, using a new
multivariate version of Stein's method of exchangeable pairs due to
S.\ Chatterjee and E.\ Meckes \cite{CM}. Our results show that, in
contrast to the situation for Dvoretzky's theorem, in the
measure-theoretic setting one can identify specific well-behaved
high-dimensional projections for large classes of convex bodies.  We
consider bodies $K$ which are $1$-unconditional, or which possess all
the symmetries of a centered regular simplex.  Another purpose of this
paper is to point out how some of the methods used here improve
quantitatively some of the results of \cite{MM}. In \cite{MM} a
symmetry hypothesis was introduced which simultaneously generalizes
$1$-unconditionality and the symmetries of a regular simplex,
described in terms of a normalized tight frame of vectors. For the
sake of transparency we have preferred to treat these special cases
independently here, although that more general setting could also be
treated with the methods of this paper.

Many of the results in this area treat marginals of probability
measures $\mu$ more general than uniform measures on convex bodies; in
particular the methods of \cite{MM} apply to completely arbitrary
probability measures which satisfy the symmetry hypotheses. One common
generalization, treated in \cite{Klartag,Klartag2} for example, is to
\emph{log-concave} measures, i.e., measures with a logarithmically
concave density with respect to Lebesgue measure. This is a natural
setting since marginals of log-concave measures are again
log-concave. While some of the methods of this paper apply to general
probability measures, we have chosen to restrict to the log-concave
case, in which stronger results are possible.

The arguments in this paper are a synthesis of the methods of the
papers \cite{MM,CM,Klartag,Klartag3}. The proofs of the main results
generalize the arguments of \cite{MM} in order to apply an abstract
normal approximation result in \cite{CM}. In order to derive stronger
results for log-concave measures, we apply a concentration result from
\cite{Klartag3} and adapt a smoothing argument from \cite{Klartag,EK}.

\medskip

The rest of this paper is organized as follows. After defining some
notation and terminology, in Section \ref{S:results} we state and
discuss our main results. Section \ref{S:smoothing} presents and
develops our tools. Finally, in Section \ref{S:proofs} we prove our
main results and make some final remarks about our methods.

\medskip

%%%%%%%%%%%%%%%%%%%%%%%%%%%%%%%%%%%%%%%%%%%%%%%%
\subsection*{Notation and terminology}

It will be convenient to frame our results in terms of random vectors
rather than probability measures.  We use $\Prob$ and $\E$ to stand
for probability and expectation respectively, and denote by $\E [Y|X]$
the conditional expectation of $Y$ given the value of $X$.

Throughout this paper $X = (X_1, \ldots, X_n)$ will be a random vector
in $\R^n$, $n\ge 2$. A random vector is \emph{$1$-unconditional} if
its distribution is invariant under reflections in the coordinate
hyperplanes of $\R^n$.  By $Z$ we denote a standard Gaussian random
vector in $\R^n$ with density
\[
\varphi_1(x) = \frac{1}{(2\pi)^{n/2}} e^{-|x|^2/2}
\]
with respect to Lebesgue measure, or a standard Gaussian random
variable in $\R$; the usage should be clear from context.  A random
vector $X$ is called \emph{isotropic} if it has mean $0$ and identity
covariance:
\[
\E X= 0, \qquad \E X_i X_j = \delta_{ij}.
\]
Observe that if $X$ is isotropic then $\E |X|^2 = n$. Isotropicity is
a natural assumption in this setting since it is preserved by
orthogonal projections and $Z$ is isotropic; see \cite{Milman} however
for recent work demonstrating that a nonisotropic affine image of $X$
is more useful in some contexts.

The total variation metric on the distributions of random vectors
in $\R^n$ may be defined by the two equivalent expressions:
\begin{equation*}\begin{split}
d_{TV} (X,Y) &= 2\sup \big\{ |\Prob[X\in A] - \Prob[Y\in A]| :
  A \subset \R^n \mbox{ measurable} \big\} \\
&= \sup \big\{ |\E f(X) - \E f(Y)|\ : \ f \in C_c(\R^n),\ 
\|f\|_\infty \le 1\big\}.
\end{split}\end{equation*}
The normalization stated here is the conventional one in analysis and
differs by a factor of $2$ from a normalization used frequently in
probability texts. Note that $d_{TV}(X,Y) = \|f-g\|_1$ if $X$ and $Y$
possess densities $f$ and $g$ respectively.  The ($L_1$-)Wasserstein
metric is defined by requiring test functions to be Lipschitz instead
of bounded:
\[
d_1 (X,Y) = \sup \big\{ |\E f(X) - \E f(Y)| \ : \forall x,y \in \R^n, \ 
  |f(x)-f(y)| \le |x-y| \big\}.
\]
Note that $d_1$ metrizes a weaker topology on probability measures
than $d_{TV}$, but quantitative results for these two metrics are not
directly comparable. In particular, $d_{TV}(X,Y) \le 2$ always, but
the typical order of magnitude of $d_1(X,Y)$ is $\sqrt{n}$.

For $x\in\R^n$, $\|x\|_p = (\sum_{i=1}^n |x_i|^p)^{1/p}$ if $1\le p
<\infty$.  Except where
noted, symbols $c, C$, etc.\ denote universal constants,
independent of $n$, $k$, and the distribution of $X$, which may differ
in value from one appearance to another.

\bigskip

%%%%%%%%%%%%%%%%%%%%%%%%%%%%%%%%%%%%%%%%
\section{Statements of the main results}\label{S:results}

Let $\theta_i = (\theta_i^1, \ldots, \theta_i^n)$, $1\le i
\le k$, be a fixed collection of $k$ orthonormal vectors. Given an
isotropic random vector $X\in\R^n$, define
\begin{equation}\label{E:W}
W_i = \inprod{X}{\theta_i}.
\end{equation}
Then $W=(W_1,\ldots,W_k)\in\R^k$ is isotropic, and is essentially the
same as $P_E(X)$, where $E$ is spanned by
$\theta_1,\ldots,\theta_k$. More concretely, $W = T(P_E(X))$, where
$T:\R^n\to \R^k$ is the partial isometry given by the matrix whose
$i$th row is $\theta_i$.  Theorems \ref{T:1-symm} and \ref{T:simplex}
give bounds on the Wasserstein and total variation distance of $W$
from a standard Gaussian random vector $Z\in\R^k$.

\begin{thm}\label{T:1-symm}
Let $X\in \R^n$ be $1$-unconditional, log-concave, and isotropic, and
let $W\in\R^k$ be as defined in \eqref{E:W}. Then
\[
d_1 (W,Z) \le 14 \sqrt{k\sum_{i=1}^k \|\theta_i\|_4^2}
\]
and
\[
d_{TV}(W,Z) \le C k^{5/6} \bigg(\sum_{i=1}^k \|\theta_i\|_4^2\bigg)^{1/3}.
\]
\end{thm}

Before stating our other main results we will makes some remarks to
put the bounds in Theorem \ref{T:1-symm} in perspective. To begin,
assume for the moment that
\begin{equation}\label{E:diagonal}
|\theta_i^\ell| = n^{-1/2} \quad \forall i,\ell.
\end{equation}
Theorem \ref{T:1-symm} then shows
\begin{equation}\label{E:rough}
d_1(W,Z) \le 14 \frac{k}{n^{1/4}}
\end{equation}
and
\begin{equation}\label{E:rough-TV}
d_{TV}(W,Z) \le C \frac{k^{7/6}}{n^{1/6}}.
\end{equation}
In particular, $d_1(W,Z) \ll 1$ as soon as $k \ll n^{1/4}$ and
$d_{TV}(W,Z) \ll 1$ as soon as $k \ll n^{1/7}$.

If $n$ is a power of $2$ then \eqref{E:diagonal} will be satisfied if
$\{\sqrt{n}\theta_i:1\le i\le k\}$ are the first $k$ vectors in a
Walsh basis for $\R^n$. For arbitrary $n$ it is not necessarily
possible to satisfy \eqref{E:diagonal} for $\theta_1,\ldots,\theta_k$
orthogonal. However, up to the values of constants, \eqref{E:rough}
and \eqref{E:rough-TV} will be satisfied by letting
$\{\sqrt{m}\theta_i:1\le i \le k\}$ be the first $k$ vectors in a
Walsh basis for $\R^m\subseteq \R^n$, where $m$ is the largest power of
$2$ not exceeding $n$ (so that $m>n/2$). In fact, a result of de
Launey \cite{deLauney} shows that one can also obtain such a so-called
partial Hadamard basis $(\theta_i)$ when $m$ is the largest multiple
of $4$ not exceeding $n$, as long as $k \le cm$ for some absolute
constant $0<c<1$. Observe that this latter condition is necessary
anyway for the bounds in \eqref{E:rough} and \eqref{E:rough-TV} to be
nontrivial.

Moreover, at the expense of the value of the constants which appear,
\eqref{E:rough} and \eqref{E:rough-TV} hold for \emph{any} orthonormal
basis $\theta_1,\ldots,\theta_k$ of \emph{most} $k$-dimensional
subspaces $E \subseteq \R^n$. This statement can be made precise using a
concentration inequality on the Grassmann manifold $G_{n,k}$ due to
Gordon \cite{Gordon}, cf.\ \cite[Lemma 16]{MM}, although we do not do
so here.

The error bounds in Theorem \ref{T:1-symm} depend on a recent optimal
concentration result for 1-unconditional log-concave random vectors
due to Klartag \cite{Klartag3}, given as Proposition
\ref{T:Klartag-var} below. Klartag used a more general version of that
result to give a sharp estimate on Gaussian approximation with respect
to Kolmogorov distance (maximum difference between distribution
functions) in the setting of Theorem \ref{T:1-symm} when $k=1$; the
typical error is of the order $1/n$. Using a smoothing lemma from
\cite{BHVV} this implies a total variation estimate (which may not be
sharp) of the order $\sqrt{\frac{\log n}{n}}$.

A Wasserstein distance estimate as in Theorem \ref{T:1-symm} can be
proved without the assumption of log-concavity, at the expense of
explicitly involving $\Var (|X|^2)$ in the bound, and (for technical
reasons) making some stronger symmetry assumption on the distribution
of $X$.  The smoothing arguments involved in proving the total
variation estimate, however, depend more crucially on log-concavity.

\medskip 

We now proceed to our other main results.

\begin{thm}\label{T:simplex}
Let $X\in\R^n$ be uniformly distributed in a regular
simplex
\begin{equation}\label{E:simplex}
\Delta_n = \sqrt{n(n+2)}\conv \{v_1,\ldots, v_n\},
\end{equation}
where $|v_i| = 1$ for $1\le i\le n+1$, and let $W\in\R^k$ be as
defined in \eqref{E:W}.  Then
\[
d_1(W,Z) \le 20 \sqrt{k \sum_{i=1}^k
\sqrt{\sum_{\ell=1}^{n+1} \inprod{\theta_i}{v_\ell}^4}}
\]
and
\[
d_{TV}(W,Z) \le C k^{5/6}\left(\sum_{i=1}^k
\sqrt{\sum_{\ell=1}^{n+1} \inprod{\theta_i}{v_\ell}^4}
\right)^{1/3}.
\]
\end{thm}

Theorem \ref{T:simplex} shows that $W$ is approximately normal as long
as the vertices of $\Delta_n$ are not close to the subspace $E = {\rm
span}(\theta_1,\ldots,\theta_k)$. By the remarks following Theorem
\ref{T:1-symm} above and \cite[Corollary 6]{MM}, Theorem
\ref{T:simplex} shows that for a typical subspace $E$, $d_{TV}(W,Z)
\le c\frac{k^{7/6}}{n^{1/6}}$.  The same proof as for Theorem
\ref{T:simplex} yields similar results for random vectors with other
distributions invariant under the symmetry group of a regular simplex.

\medskip

Our last main result improves the typical dependence on $n$ of the
total variation bound of Theorem \ref{T:simplex} in the case that $k=1$.

\begin{thm}\label{T:univ-sncp}
Let $X$ be uniformly distributed in a regular simplex $\Delta_n$
as defined in \eqref{E:simplex}, let $\theta\in S^{n-1}$ be fixed,
and let $W=\inprod{X}{\theta}$. Then
\[ d_{TV}(W,Z) \le C \sqrt{\sum_{i=1}^{n+1} |\inprod{\theta}{v_i}|^3},\]
where $Z$ is a standard Gaussian random variable.
\end{thm}

For a typical $\theta\in S^{n-1}$, we obtain here $d_{TV}(W,Z) \le c
n^{-1/4}$. This also improves an error bound given in \cite{MM}; see
the remarks at the end of this paper for further details.

\bigskip

%%%%%%%%%%%%%%%%%%%%%%%%%%%%%%%%%%%%%%%%
\section{Smoothing and abstract Gaussian approximation 
theorems}\label{S:smoothing}

For $n\ge 1$ and $t>0$, define $\varphi_t:\R^n \to \R$ by
\[ \varphi_t(x) = \frac{1}{(2\pi t^2)^{n/2}} e^{-|x|^2/2t^2},\]
so $\varphi_t$ is the density of $tZ$, where $Z\in\R^n$ is a standard
Gaussian random vector. A well-known consequence of the
Pr\'ekopa-Leindler inequality \cite{Leindler,Prekopa} (or see
\cite{Gardner}) is that the convolution of integrable log-concave
functions is log-concave; hence in particular $f*\varphi_t$ is
log-concave for any log-concave probability density $f:\R^n\to\R_+$.
Furthermore it is well-known that
\begin{equation*}
\|f*\varphi_t - f\|_1 \to 0 \quad \mbox{as $t\to 0$}
\end{equation*}
for any integrable $f$. Thus log-concave random variables are
arbitrarily well approximated, in the total variation metric, by
log-concave random vectors with smooth densities. The statements of
Section \ref{S:results} involving the total variation metric rely on a
quantitative version of this observation. We say that $f:\R^n\to \R_+$
is isotropic if it is the density of an isotropic random vector in
$\R^n$.  The following is a sharp version of Lemma 5.1 in
\cite{Klartag}.

\begin{prop}\label{T:smoothing}
Let $f:\R^n\to \R_+$ an isotropic and log-concave probability density.
\begin{enumerate}
\item \label{I:smooth-1} If $n=1$ then
\[ \| f * \varphi_t - f \|_1 \le 2\sqrt{2}\ t\]
for all $t\ge 0$.
\item (Klartag-Eldan) \label{I:smooth-2} If $n\ge 2$ then
\[ 
\| f * \varphi_t - f \|_1 \le cnt 
\]
for all $t\ge 0$, where $c>0$ is an absolute constant.
\end{enumerate}
\end{prop}

Proposition \ref{T:smoothing}(\ref{I:smooth-2}) was conjectured in an
earlier version of this paper which also proved a weaker version of
that estimate by optimizing over some of the parameters in the proof
of \cite[Lemma 5.1]{Klartag}.  After that version of this paper was
posted on arxiv.org, Klartag proved the conjecture; the proof appears
in \cite[Section 5]{EK}.

For arbitrary $f$, Proposition \ref{T:smoothing} is sharp up to the
values of the constants $2\sqrt{2}$ and $c$. For the particular case
$f=\varphi_1$, one can show
\begin{equation}\label{E:heat-normal}
\|\varphi_1 * \varphi_t - \varphi_1\|_1 
	= \|\varphi_{\sqrt{1+t^2}} - \varphi_1\|_1 < \sqrt{2n}\ t
\end{equation}
for all $t\ge 0$ and any $n$ (cf.\ \cite[Lemma 4.9]{Klartag}, or the
proof of Proposition \ref{T:smoothing}(\ref{I:smooth-1}) below).

\begin{proof}[Proof of Proposition \ref{T:smoothing}(\ref{I:smooth-1})]
First, we can assume that $f$ is smooth and everywhere positive, for
example by convolving $f$ with $\varphi_\varepsilon$, rescaling for
isotropicity, and letting $\varepsilon \to 0$. A special case of a
result of Ledoux \cite[formula (5.5)]{Ledoux3} about the heat
semigroup on a Riemannian manifold implies that
\begin{equation}\label{E:heat}
\| f * \varphi_t - f\|_1 \le \sqrt{2}\ t\ \| f' \|_1
\end{equation}
for any $t\ge 0$. Since $f$ is log-concave it is unimodal, i.e., there
exists an $a\in \R$ such that $f'(s) \ge 0$ for $s\le a$ and $f'(s)
\le 0$ for $s\ge a$, and since $f$ is also isotropic we have $f(a) \le
1$ (see e.g.\ \cite[Lemma 5.5(a)]{LV}). Therefore
\[ \|f'\|_1 = \int_{-\infty}^a f'(s) \ ds
		- \int_a^\infty f'(s)\ ds = 2f(a) \le 2,\]
which proves the first claim.
\end{proof}

\bigskip

The following abstract Gaussian approximation theorem was proved by
Stein in \cite{Stein}; the version stated here incorporates a slight
improvement in the constants proved in \cite{BRS}. Recall that a pair
of random variables $(W,W')$ is called exchangeable if the joint
distribution of $(W,W')$ is the same as the distribution of $(W',W)$.

\begin{prop}[Stein]\label{T:Stein}
Suppose that $(W,W')$ is an exchangeable pair of random variables such
that $\E W = 0$, $\E W^2 = 1$, and $\E[W'-W|W] = -\lambda W$.  Then
\[
|\E g(W) - \E g(Z)| 
	\le \frac{\|g\|_\infty}{\lambda}\sqrt{\Var \E [(W'-W)^2 | W]}
	+ \frac{\|g'\|_\infty}{4\lambda} \E |W'-W|^3
\]
for any $g\in C_c^\infty(\R)$, where $Z\in\R$ denotes a standard
Gaussian random variable.
\end{prop}

Stein used a smoothing argument to derive a version of Proposition
\ref{T:Stein} for the Kolmogorov distance, which was the main tool in
the proofs of most of the results of \cite{MM}. Estimates for total
variation distance for log-concave distributions were obtained in
\cite{MM} by combining the Kolmogorov distance estimates with
\cite[Theorem 3.3]{BHVV}, which entails an additional loss in the
error bound. Here we use Proposition \ref{T:smoothing}(\ref{I:smooth-1})
to obtain a version of Proposition \ref{T:Stein} for total variation
distance and log-concave distributions, which matches Stein's bound
for Kolmogorov distance used in \cite{MM}; this is the main technical
tool in the proof of Theorem \ref{T:univ-sncp}.

\begin{cor}\label{T:TV-univ}
Suppose, in addition to the hypotheses of Proposition \ref{T:Stein},
that $W$ is log-concave. Then
\[
d_{TV}(W,Z) \le
	\frac{1}{\lambda}\sqrt{\Var \E [(W'-W)^2|W]}
	+ 2 \sqrt{\frac{1}{\lambda} \E |W'-W|^3}.
\]
\end{cor}

It follows from the proof of Proposition
\ref{T:smoothing}(\ref{I:smooth-1}) that Corollary \ref{T:TV-univ}
only requires $W$ to have a bounded unimodal density with respect to
Lebesgue measure; the coefficient $2$ in the r.h.s.\ should be
replaced by a constant depending on the maximum value of the density.

\begin{proof}
Let $g \in C_c(\R)$ with $\|g\|_\infty \le 1$, and let $f$ denote
the density of $W$. Assume for now that $f$ is smooth. Given
$t>0$, define $h=g*\varphi_t$. To begin, observe that
\[ \|h\|_\infty \le \|g\|_\infty \|\varphi_t\|_1 \le 1, \]
\[ \|h'\|_\infty = \| g*(\varphi_t')\|_\infty
	\le \|g\|_\infty \|\varphi_t'\|_1 \le 
	\frac{1}{t}\sqrt{\frac{2}{\pi}}\]
Proposition \ref{T:Stein} applied to $h$ implies that
\begin{equation}\label{E:TV-univ-1}
|\E h(W) - \E h(Z)| \le \frac{1}{\lambda}\sqrt{\Var\E[(W'-W)^2|W]}
	+ \frac{1}{2\sqrt{2\pi}\ \lambda t} \E |W'-W|^3.
\end{equation}

Next, by Proposition \ref{T:smoothing}(\ref{I:smooth-1}),
\begin{equation}\label{E:TV-univ-2}\begin{split}
|\E h(W) - \E g(W)| &= \left|\int [g*\varphi_t(s)-g(s)] f(s)\ ds
	\right|
= \left| \int_{\R} [f*\varphi_t(s)-f(s)]g(s)\ ds \right| \\
&\le \|g\|_\infty \| f*\varphi_t - f\|_1
\le 2\sqrt{2}\ t. 
\end{split}\end{equation}
Similarly, by \eqref{E:heat}
\begin{equation}\label{E:TV-univ-3}
|\E h(Z) - \E g(Z)| 
\le \|g\|_\infty \| \varphi_1*\varphi_t - \varphi_1\|_1
< \sqrt{2} \|\varphi_1'\|_1 \ t \ \le  \frac{2}{\sqrt{\pi}}\ t.
\end{equation}
The only reason for using \eqref{E:heat} directly instead of 
applying Proposition \ref{T:smoothing}(\ref{I:smooth-1}) here is to obtain a
slightly better constant. Combining \eqref{E:TV-univ-1},
\eqref{E:TV-univ-2}, and \eqref{E:TV-univ-3} yields
\begin{equation}\label{E:TV-univ-4}
|\E g(W) - \E g(Z)| \le \frac{1}{\lambda}\sqrt{\Var\E[(W'-W)^2|W]}
	+ \frac{1}{2\sqrt{2\pi}\ \lambda t} \E |W'-W|^3
	+ 2\left(\sqrt{2}+\frac{1}{\sqrt{\pi}}\right) t.
\end{equation}
The corollary, under the assumption that $f$ is smooth, now follows by
optimizing in $t$. The coefficient of $2$ given in the second term in
the statement of the corollary is not optimal and is given as such for
simplicity.

The corollary can be reduced to the smooth case with a convolution
argument as for Proposition \ref{T:smoothing}, although it is slightly
more complicated because it is necessary to smooth not only $f$ but
the exchangeable pair $(W,W')$. To do this, let $Z_1, Z_2$ be
standard Gaussian random variables independent of each other and
of $(W,W')$, set $Z=Z_1$, and set
\[
Z' = (1-\lambda)Z_1 + \sqrt{2\lambda - \lambda^2} \ Z_2.
\]
Then $(Z,Z')$ is an exchangeable pair and $\E[Z'-Z|Z] = -\lambda Z$.
Now for $\varepsilon>0$ let
\[
W_\varepsilon = \frac{1}{\sqrt{1+\varepsilon^2}}(W+\varepsilon Z),
\quad
W'_\varepsilon = \frac{1}{\sqrt{1+\varepsilon^2}}(W'+\varepsilon Z').
\]
Then $(W_\varepsilon,W'_\varepsilon)$ is an exchangeable pair that
satisfies all the hypotheses of the corollary (log-concavity follows
from the Pr{\'e}kopa-Leindler inequality), and $W_\varepsilon$ has a
smooth density. Applying the corollary to $(W_\varepsilon,
W'_\varepsilon)$ and letting $\varepsilon\to 0$ yields the general
case.
\end{proof}

\medskip

The main technical tool in the proofs of Theorems
\ref{T:1-symm} and \ref{T:simplex} is the following multivariate version
of Proposition \ref{T:Stein}, recently proved by S.\ Chatterjee and
E.\ Meckes in \cite{CM}. For a smooth function $f:\R^k
\to \R$, we denote by $M_1(f) = \| |\nabla f| \|_\infty$ the Lipschitz
constant of $f$ and 
\[
M_2(f) = \| \| \nabla^2 f \| \|_\infty
\]
the maximum value of the operator norm of the Hessian of $f$, or
equivalently the Lipschitz constant of $\nabla f:\R^k \to \R^k$.

\begin{prop}[Chatterjee and E.\ Meckes]\label{T:CM}
Let $W$ and $W'$ be identically distributed random vectors in $\R^k$
defined on a common probability space.  Suppose that for some constant
$\lambda > 0$ and random variables $E_{ij}$, $1\le i,j \le k$,
\begin{align*}
\E \big[W'_i - W_i \big| W \big] &= - \lambda W_i, \\
\E \big[(W'_i - W_i) (W'_j - W_j) \big| W \big] 
  &= 2\lambda \delta_{ij} + E_{ij}.
\end{align*}
Then
\begin{equation*}
\big| \E f(W) - \E f(Z) \big|
\le \frac{M_1(f)}{\lambda}
  \E \sqrt{\sum_{i,j=1}^k E_{ij}^2}
  + \frac{\sqrt{2\pi}M_2(f)}{24 \lambda} \E |W'-W|^3
\end{equation*}
for any smooth $f:\R^k \to \R$, where $Z\in\R^k$ is
a standard Gaussian random vector.
\end{prop}

Note that the normalization for $E_{ij}$ used here differs from that
in the statement of \cite[Theorem 4]{CM} by a factor of $\lambda$.
An earlier version of this paper was posted to arxiv.org which was
based on an earlier version of Proposition \ref{T:CM}. The version
given above allows improved estimates in Theorem \ref{T:1-symm}.

Convolution arguments similar to those in the proof of Corollary
\ref{T:TV-univ} yield bounds on Wasserstein and total variation
distances.

\begin{cor}\label{T:Wass-TV}
Under the same hypotheses as Proposition \ref{T:CM},
\[
d_1(W,Z) \le \frac{1}{\lambda} \E \sqrt{\sum_{i,j=1}^k E_{ij}^2}
  + k^{1/4}
  \sqrt{\frac{2}{3\lambda} \E |W' - W|^3}.
\]
If moreover $W$ is log-concave, then
\[
d_{TV}(W,Z) \le C
  \left( \frac{k}{\lambda} \E \sqrt{\sum_{i,j=1}^k E_{ij}^2}
    +\frac{k^2}{\lambda} \E |W'-W|^3\right)^{1/3},
\]
where $C>0$ is an absolute constant.
\end{cor}

\begin{proof}
To prove the first claim, let $g:\R^k \to \R$ be $1$-Lipschitz, and
define $h=g*\varphi_t$ for $t>0$. Standard calculations show
\begin{equation}\label{E:smooth-L}
M_1(h) \le M_1(g) \|\varphi_t\|_1 \le 1 
\end{equation}
and
\begin{equation}\label{E:smooth-2}
M_2(h) \le M_1(g) \sup_{\theta\in S^{n-1}} 
  \|\inprod{\nabla \varphi_t}{\theta}\|_1
  \le \sqrt{\frac{2}{\pi}} \frac{1}{t}.
\end{equation}
Note that $\E h(W) = \E g(W+tZ)$, where $Z\in\R^k$ is a standard
Gaussian random vector independent of $W$, which implies 
\begin{equation}\label{E:smooth-W}
\big|\E h(W) - \E g(W) \big| \le \E |tZ| \le t\sqrt{k}
\end{equation}
since $g$ is $1$-Lipschitz, and similarly
\begin{equation}\label{E:smooth-Z}
\big|\E h(Z)-\E g(Z)\big| \le t\sqrt{k}.
\end{equation}
The claim follows by applying Proposition \ref{T:CM} to $h$, using
\eqref{E:smooth-L}, \eqref{E:smooth-2}, \eqref{E:smooth-W}, and
\eqref{E:smooth-Z}, and optimizing in $t$.

\medskip

The proof of the second claim is similar.
As in the proof of Corollary \ref{T:TV-univ} we may assume that
$W$ has a smooth density $f$. Let $g\in
C_c(\R^k)$ with $\| g \|_\infty \le 1$, and again define $h=g*\varphi_t$ for
$t>0$. By standard calculations,
\begin{equation}\label{E:smoother-L}
M_1(h) \le \|g\|_\infty \sup_{\theta\in S^{n-1}}
  \| \inprod{\nabla \varphi_t}{\theta}\|_1
  \le \sqrt{\frac{2}{\pi}} \frac{1}{t}
\end{equation}
and
\begin{equation}\label{E:smoother-2}
M_2(h) \le \|g\|_\infty \sup_{\theta\in S^{n-1}}
  \| \inprod{(\nabla^2 \varphi_t)\theta}{\theta} \|_1
  \le \frac{\sqrt{2}}{t^2}.
\end{equation}
Proposition \ref{T:smoothing}(\ref{I:smooth-2}) and the identity
\[ \E h(W) = \int g*\varphi_t (x) \ f(x)\ dx
  = \int f*\varphi_t (x) \ g(x) \ dx\]
imply that
\begin{equation}\label{E:smoother-W}
\big|\E h(W) - \E g(W) \big| \le \|g\|_\infty 
  \| f*\varphi_t - f \|_1 \le ckt.
\end{equation}
Similarly, 
\begin{equation}\label{E:smoother-Z}
\big|\E h(Z) - \E g(Z) \big| \le \|g\|_\infty 
  \|\varphi_1*\varphi_t - \varphi_1 \|_1 \le ckt
\end{equation}
This last estimate can be improved to $\sqrt{2k}t$ using
\eqref{E:heat-normal}, but there is no advantage to doing so here.

Applying Proposition \ref{T:CM} to $h$ and using \eqref{E:smoother-L},
\eqref{E:smoother-2}, \eqref{E:smoother-W}, and \eqref{E:smoother-Z}
yields
\begin{equation}\label{E:smoother-TV1}
d_{TV}(W,Z) = \sup_{g\in C_c(\R^k),\ \|g\|_\infty\le 1}
\big|\E g(W) - \E g(Z) \big|
  \le \frac{A}{t^2} + \frac{B}{t} + ckt
\end{equation}
for any $t>0$, where
\[
A = \frac{\sqrt{\pi}}{12\lambda} \E |W' - W|^3, \qquad
B = \frac{\sqrt{2}}{\pi}\frac{1}{\lambda} \E \sqrt{\sum_{i,j=1}^k E_{ij}^2}.
\]
Although it is not straightforward to optimize the r.h.s.\ of
\eqref{E:smoother-TV1} precisely, this is simplified by noting that
$d_{TV}(W,Z) \le 2$ always. Therefore \eqref{E:smoother-TV1}
is vacuously true for $t \ge 2/(ck)$, and so
\[
d_{TV}(W,Z) \le \frac{A+2B/(ck)}{t^2} + ckt
\]
for any $t>0$. Optimizing this latter expression in $t$ yields
\[
d_{TV}(W,Z) \le C \big(k^2 A + k B\big)^{1/3},
\]
from which the result follows.
\end{proof}

\bigskip 

%%%%%%%%%%%%%%%%%%%%%%%%%%%%%%%%%%%%%%%%
\section{Proofs of the main results}\label{S:proofs}

In this section we prove Theorems \ref{T:1-symm} and \ref{T:simplex}
from Corollary \ref{T:Wass-TV}, and indicate how Theorems
\ref{T:univ-sncp} may be proved from Corollary \ref{T:TV-univ}.  The
arguments mostly generalize the proofs of \cite{MM}.

First, observe that the Pr\'ekopa-Leindler theorem
\cite{Prekopa,Leindler} implies that marginals of log-concave measures
are log-concave. Therefore when $X$ is log-concave, $W$ is log-concave
as well, and the second estimate of Corollary \ref{T:Wass-TV} may be
applied. This fact will be used without further comment in all the
proofs in this section.

Second, we state a version of Klartag's concentration result for
unconditional convex bodies.

\begin{prop}[Klartag]\label{T:Klartag-var}
If $X$ is isotropic, unconditional, and log-concave, and
$a_1,\dotsc,a_n \in \R$, then
\begin{equation*}
\Var\left(\sum_{\ell=1}^n a_\ell X_\ell^2\right) 
  \le 32 \sum_{\ell=1}^n a_\ell^2.
\end{equation*}
\end{prop}

Proposition \ref{T:Klartag-var} is essentially a special case of
\cite[Lemma 4]{Klartag3}, which is stated with the additional
assumption that $a_1,\dotsc,a_n \ge 0$. For the precise constants
which appear here see the comments following the proof in
\cite{Klartag3}; an extra factor of $2$ is introduced to allow
negative coefficients by observing that
\[
\Var(X+Y) \le 2(\Var X + \Var Y)
\]
for any pair of random variables $X$ and $Y$.

We now proceed with the proofs of our main results.

\begin{proof}[Proof of Theorem \ref{T:1-symm}]
To construct $W'$ appropriately coupled with $W$, we first define
$X'$ by reflecting $X$ in a randomly chosen coordinate hyperplane,
and then let $W'_i = \inprod{X'}{\theta_i}$.
By the $1$-unconditionality of $X$, $X$ and $X'$ are identically
distributed and hence so are $W$ and $W'$.

More precisely, let $I$ be a random variable chosen uniformly from
$\{1,\ldots,n\}$ and independently from the random vector $X$. Then
\[ X' = X - 2X_I e_I,\]
where $e_i$ is the $i$th standard basis vector in $\R^n$, and
\[ W'_i = \inprod{X-2X_I e_I}{\theta_i} = W_i - 2 \theta_i^I X_I.\]
It follows that
\[
\E\big[ W'_i - W_i \big| X \big] = -\frac{2}{n}\sum_{\ell = 1}^n
  \theta_i^\ell X_\ell = -\frac{2}{n}W_i
\]
and
\[
\E\big[(W'_i-W_i)(W'_j-W_j)\big| X \big]
= \frac{4}{n}\sum_{\ell=1}^n \theta_i^\ell \theta_j^\ell X_\ell^2.
\]
Therefore we may apply Corollary \ref{T:Wass-TV} with $\lambda =
\frac{2}{n}$ and
\[
E_{ij} = \frac{4}{n}\left(\sum_{\ell=1}^n \theta_i^\ell \theta_j^\ell
  X_\ell^2 - \delta_{ij} \right).
\]

Now by Jensen's inequality, Proposition \ref{T:Klartag-var}, and the
Cauchy-Schwarz inequality, 
\begin{equation*}\begin{split}
\frac{1}{\lambda}\E \sqrt{\sum_{i,j=1}^k E_{ij}^2}
&\le \frac{1}{\lambda}\sqrt{\sum_{i,j=1}^k \E E_{ij}^2}
= 2\sqrt{\sum_{i,j=1}^k \Var \left(\sum_{\ell=1}^n
  \theta_i^\ell \theta_j^\ell X_\ell^2\right)}\\
& \le 8 \sqrt{2 \sum_{i,j=1}^k 
  \sum_{\ell=1}^n (\theta_i^\ell)^2(\theta_j^\ell)^2}
\le 8 \sqrt{2} \sum_{i=1}^k \| \theta_i\|_4^2.
\end{split}\end{equation*}

By the triangle inequality for the $L_{3/2}$ norm and a precise
version of Borell's lemma (found, e.g., in \cite{MP}),
\begin{equation}\begin{split} \label{E:3rdMoment}
\E |W - W|^3 &= \E \bigg( \sum_{i=1}^k |W_i'-W_i|^2 \bigg)^{3/2}
  \le \bigg(\sum_{i=1}^k \big(\E|W_i'-W_i|^3 \big)^{2/3} \bigg)^{3/2} \\
&=\frac{8}{n} \left(\sum_{i=1}^k \bigg(\sum_{\ell=1}^n
  |\theta_i^\ell|^3 \E |X_\ell|^3 \bigg)^{2/3} \right)^{3/2}
\le \frac{12\sqrt{2}}{n}\bigg(\sum_{i=1}^k \|\theta_i\|_3^2 \bigg)^{3/2}.
\end{split}\end{equation}

By the standard estimates between $\ell_p^k$ norms and the fact that
$\|\theta_i\|_3^3 \le |\theta_i| \|\theta_i\|_4^2 = \|\theta\|_4^2$,
\begin{equation}\label{E:simplifying}
\bigg(\sum_{i=1}^k \|\theta_i\|_3^2 \bigg)^{3/2}
  \le \sqrt{k} \sum_{i=1}^k \|\theta_i\|_3^3
  \le \sqrt{k} \sum_{i=1}^k \|\theta_i\|_4^2.
\end{equation}

Proposition \ref{T:Wass-TV} now implies the stated bound for $d_{TV}$
immediately. For the bound on $d_1$ observe also that
$\|\theta_i\|_4^2 \le |\theta_i|^2 =1$, and so $\sum_{i=1}^k
\|\theta_i\|_4^2 \le k$.
\end{proof}

\medskip

By using Proposition \ref{T:CM} directly, the proof of Theorem
\ref{T:1-symm} above yields better bounds on the distance
\[ 
d_2(W,Z) = \sup\big\{ |\E f(W) - \E f(Z)| : M_1(f), M_2(f)\le 1
\big\}.
\]
In particular, as in the remarks following the statement of Theorem
\ref{T:1-symm}, under the conditions of that theorem, typical
$k$-dimensional marginals are nearly Gaussian with respect to $d_2$ if
$k \ll n^{1/3}$.  The same remark applies to the proof of Theorem
\ref{T:simplex} below.  While we have preferred here to work with the
more classical Wasserstein and total variation metrics, metrics like
$d_2$ based on smooth test functions are commonly used in quantifying
multivariate Gaussian approximation.

It is also worth pointing out here that \cite{MM} does prove
multivariate Gaussian approximation results, but with respect to a
weak metric referred to as $T$-distance which captures only the
behavior of $1$-dimensional marginals. Using $T$-distance yields
misleadingly good results in terms of how large $k$ may be for an
approximately Gaussian marginal, cf.\ the remarks at the very end of
\cite{Klartag}. Metrics like $d_{TV}$, $d_1$, and $d_2$ based on
regular test functions better capture high-dimensional behavior.

\medskip

\begin{proof}[Proof of Theorem \ref{T:simplex}]
In this case $X'$ is obtained by reflecting $X$ in a hyperplane spanned 
by $(n-1)$ vertices of $\Delta_n$; alternatively one may
think of this operation as transposing two vertices.

We will need the well-known facts about vertices of centered regular
simplices (which may be seen e.g.\ as consequences of John's
theorem, cf.\ \cite{Ball}) that
\begin{equation}\label{E:simplex-0}
\sum_{i=1}^{n+1} v_i = 0
\end{equation}
and
\begin{equation}\label{E:simplex-v}
\sum_{i=1}^{n+1} \inprod{x}{v_i} v_i = \frac{n+1}{n} x
\end{equation}
for any $x\in \R^n$. It will be convenient to use the notation
\[ u_{ij} = \sqrt{\frac{n}{2(n+1)}} (v_i-v_j),
\quad 1\le i,j \le n+1,\]
and $x^{ij}=\inprod{x}{u_{ij}}$ for $x\in \R^n$.
It follows from \eqref{E:simplex-0} and \eqref{E:simplex-v} that 
$|u_{ij}| = 1$ for $i\neq j$ and
\begin{equation}\label{E:simplex-u}
\sum_{\ell\neq m} x^{\ell m}u_{\ell m} = (n+1) x
\quad \forall x\in \R^n.
\end{equation}

To define $W'$ precisely, first pick a pair
$(I,J)$ of distinct elements of $\{1, \ldots, n+1\}$ uniformly and
independently of $X$. Let
\[ X' = X - 2 X^{IJ} u_{IJ}.\]
and $W'_i = \inprod{X'}{\theta_i}$ as before.
Using \eqref{E:simplex-u}, one obtains
\begin{align*}
\E \big[W'_i - W_i \big| X \big]
&= -\frac{2}{n}\inprod{X}{\theta_i} = -\frac{2}{n}W_i,\\
\E \big[(W'_i - W_i) (W'_j - W_j) \big| X \big] 
&= \frac{4}{n(n+1)} \sum_{\ell,m=1}^{n+1} \theta_i^{\ell m}
  \theta_j^{\ell m} (X^{\ell m})^2,
\end{align*}
and so Corollary \ref{T:Wass-TV} applies
with $\lambda=\frac{2}{n}$ and
\[
E_{ij} = \frac{4}{n} \left(\frac{1}{n+1} \sum_{\ell,m=1}^{n+1} 
  \theta_i^{\ell m} \theta_j^{\ell m} 
  (X^{\ell m})^2 - \delta_{ij}. \right)
\]
The relevant moments were calculated in \cite{MM}: for $\ell \neq m$
and $p \neq q$,
\begin{equation}\label{E:simplex-moments}
\E (X^{\ell m})^2 (X^{p q})^2 
= \frac{(n+1)(n+2)}{(n+3)(n+4)} 
\begin{cases}
1 & \mbox{ if } \{\ell, m\} \cap \{p, q\} = \emptyset, \\
3 & \mbox{ if } |\{\ell, m\} \cap \{p, q\}| = 1, \\
6 & \mbox{ if } |\{\ell, m\} \cap \{p, q\}| = 2,
\end{cases}
\end{equation}
and
\begin{equation} \label{E:simplex-3}
\E |X^{\ell m}|^3 < 3\sqrt{2}.
\end{equation}

In order to estimate $\E E_{ij}^2$, decompose the resulting sum of
terms involving $\E (X^{\ell m})^2 (X^{pq})^2$ according to the size
of $\{\ell, m\}\cap \{p,q\}$ and use \eqref{E:simplex-moments}:
\begin{multline}\label{E:simplex-bigsum}
\frac{(n+3)(n+4)}{(n+1)(n+2)}\sum_{\ell, m, p, q}
  \theta_i^{\ell m}\theta_i^{pq}
  \theta_j^{\ell m}\theta_j^{pq}
  \E (X^{\ell m})^2 (X^{pq})^2 \\
= \left(\sum_{\ell,m}
  \theta_i^{\ell m}\theta_j^{\ell m}\right)^2
+ 2\sum_{|\{\ell,m\}\cap\{p,q\}|=1}
  \theta_i^{\ell m}\theta_i^{pq}
  \theta_j^{\ell m}\theta_j^{pq}
+ 10\sum_{\ell, m}
  (\theta_i^{\ell m})^2(\theta_j^{\ell m})^2.
\end{multline}
In all of the sums in \eqref{E:simplex-bigsum} the indices range from
$1$ to $n+1$, with $\ell \neq m$ and $p\neq q$ and in the last term we
have also used that $u_{\ell m} = - u_{m\ell}$.  It follows from
\eqref{E:simplex-u} that
\begin{equation}\label{E:simplex-delta}
\sum_{\ell\neq m}
  \theta_i^{\ell m}\theta_j^{\ell m}
= (n+1)\inprod{\theta_i}{\theta_j} = (n+1) \delta_{ij}.
\end{equation}
By the Cauchy-Schwarz inequality, the definition of $u_{\ell m}$, and
the $\ell_4$ triangle inequality,
\begin{equation}\label{E:simplex-CS1}
\sum_{\ell \neq m} (\theta_i^{\ell m})^2 (\theta_j^{\ell m})^2
\le \sqrt{\sum_{\ell \neq m} (\theta_i^{\ell m})^4}
\sqrt{\sum_{\ell \neq m} (\theta_j^{\ell m})^4}
\le 4n 
  \sqrt{\sum_{\ell=1}^{n+1} \inprod{\theta_i}{v_\ell}^4}
  \sqrt{\sum_{\ell=1}^{n+1} \inprod{\theta_j}{v_\ell}^4}.
\end{equation}
Also by the Cauchy-Schwarz inequality,
\begin{equation}\label{E:simplex-CS2}
\sum_{\substack{|\{\ell,m\}\cap\{p,q\}|=1 \\ \ell \neq m,\ p \neq q}}
  \theta_i^{\ell m}\theta_i^{pq}
  \theta_j^{\ell m}\theta_j^{pq}
\le 
\sum_{\substack{|\{\ell,m\}\cap\{p,q\}|=1 \\ \ell \neq m,\ p \neq q}}
	(\theta_i^{\ell m})^2 (\theta_j^{\ell m})^2
\le 4n
\sum_{\ell \neq m} (\theta_i^{\ell m})^2 (\theta_j^{\ell m})^2.
\end{equation}

Combining \eqref{E:simplex-bigsum}, \eqref{E:simplex-delta},
\eqref{E:simplex-CS1}, and \eqref{E:simplex-CS2},
\begin{equation*}
\begin{split}
\E E_{ij}^2
&= \frac{16}{n^2} \Bigg( \frac{1}{(n+1)^2}
\sum_{\ell, m, p, q}
  \theta_i^{\ell m}\theta_i^{pq}
  \theta_j^{\ell m}\theta_j^{pq}
  \E (X^{\ell m})^2 (X^{pq})^2 - \delta_{ij} \Bigg) \\
&\le \frac{512}{n^2}  \sqrt{\sum_{\ell=1}^{n+1} 
	\inprod{\theta_i}{v_\ell}^4}
  	\sqrt{\sum_{\ell=1}^{n+1} \inprod{\theta_j}{v_\ell}^4},
\end{split}\end{equation*}
and therefore the first error term in Corollary \ref{T:Wass-TV} is
bounded by
\[
\frac{1}{\lambda}\E\sqrt{\sum_{i,j=1}^k E_{ij}^2} \le 
8\sqrt{2} \sum_{i=1}^k \sqrt{\sum_{\ell=1}^{n+1}
	\inprod{\theta_i}{v_\ell}^4}.
\]

To bound the second error term, we begin as in \eqref{E:3rdMoment},
using \eqref{E:simplex-3}, the definition of $u_{\ell m}$, and the
$\ell_3$ triangle inequality to obtain
\begin{equation*}\begin{split}
\E |W - W|^3 & \le 
  \bigg(\sum_{i=1}^k \big(\E|W_i'-W_i|^3 \big)^{2/3} \bigg)^{3/2}
=\frac{8}{n(n+1)} \left(\sum_{i=1}^k \bigg(\sum_{\ell\neq m}
  |\theta_i^{\ell m}|^3 \E |X^{\ell m}|^3 \bigg)^{2/3} \right)^{3/2}\\
&\le \frac{12\sqrt{n}}{(n+1)^{5/2}} \left(\sum_{i=1}^k \bigg(
  \sum_{\ell \neq m}|\inprod{\theta_i}{v_\ell} - \inprod{\theta_i}{v_m}|^3
  \bigg)^{2/3} \right)^{3/2}\\
&\le \frac{96\sqrt{n}}{(n+1)^{3/2}} \left(\sum_{i=1}^k 
  \bigg(\sum_{\ell=1}^{n+1} |\inprod{\theta_i}{v_\ell}|^3\bigg)^{2/3}
  \right)^{3/2}.
\end{split}\end{equation*}
The error bounds from Proposition \ref{T:Wass-TV} are now simplified
similarly as in the proof of Theorem \ref{T:1-symm}. For each $i$
define $x_i = (x_i^1,\ldots, x_i^{n+1})$ by $x_i^\ell =
\inprod{\theta_i}{v_\ell}$. By \eqref{E:simplex-v},
\[
|x_i|^2 = \sum_{\ell=1}^{n+1} \inprod{\theta_i}{v_\ell}^2
  = \frac{n+1}{n} |\theta_i|^2 = \frac{n+1}{n}.
\]
Therefore $\|x_i\|_3^3 \le |x_i| \|x_i\|_4^2 \le \sqrt{\frac{n+1}{n}}
\|x_i\|_4^2$, and so by the same reasoning as in
\eqref{E:simplifying},
\[
\E |W'-W|^3 
  \le \frac{96\sqrt{k}}{n+1}\sum_{i=1}^k\sqrt{\sum_{\ell=1}^{n+1}
  \inprod{\theta_i}{v_\ell}^4};
\]
finally observe also that $\|x_i\|_4^2 \le |x|^2$ to simplify the
bound on $d_1$.
\end{proof}

\medskip

Theorem \ref{T:univ-sncp} may be proved by following the proof of
Theorem \ref{T:simplex} in the case $k=1$, applying Corollary
\ref{T:TV-univ} in place of Corollary \ref{T:Wass-TV}. Alternatively,
one can follow the proof of \cite[Corollary 6]{MM}, using Corollary
\ref{T:TV-univ} in place of the Stein's Kolmogorov distance version of
Proposition \ref{T:Stein}; this amounts to the same thing.

In \cite{MM} Stein's Kolmogorov distance version of Proposition
\ref{T:Stein} was applied for arbitrary isotropic $X$ (under various
symmetry assumptions). Total variation estimates for the log-concave
case were then deduced using \cite[Theorem 3.3]{BHVV}, which allows
Gaussian approximation estimates for log-concave random variables to
be transferred from Kolmogorov distance to total variation distance.
The present approach entails less loss in the final total variation
bound since it uses only one smoothing argument instead of two. In
general, the second approach described above to prove Theorem
\ref{T:univ-sncp} can be used to deduce total variation bounds of the
same order as the Kolmogorov distance bounds in most of the results of
\cite{MM} for log-concave random vectors.  In particular, this applies
to Theorem 1, Corollary 4(2), and parts of Corollary 5 of \cite{MM}.

In a similar fashion, using Proposition \ref{T:Stein} directly yields
versions of many of the results of \cite{MM} for the bounded Lipschitz
metric
\[
d_{BL}(X,Y) = \sup \big\{|\E f(X) -\E f(Y)| : \|f\|_\infty, |f|_L \le 1
\big\}.
\]
In general, the dominant error term for the results of \cite{MM} using
$d_{BL}$ is typically of the order $n^{-1/2}$, as opposed to the order
$n^{-1/4}$ for the Kolmogorov distance in most of the results of
\cite{MM} and for $d_{TV}$ in the present Theorem \ref{T:univ-sncp}.

In \cite{Klartag}, Klartag used another approach, based on an
application of the classical Berry-Esseen theorem, to prove a
univariate estimate in the setting of Theorem \ref{T:1-symm}.  Since
Stein's method can be used to prove the Berry-Esseen theorem, the
approach taken here and in \cite{MM} is arguably more direct, and the
total variation bounds which can be derived in this way are better
than those derived by the method of \cite{Klartag}.  However,
since the original version of this paper was written, in
\cite{Klartag3} Klartag has given a proof of an optimal result for
Kolmogorov distance, which, as discussed after the statement of
Theorem \ref{T:1-symm} above, implies sharper total variation bounds
(when $k=1$) than the methods used here. Klartag's proof is based
partly on the optimal concentration result proved in \cite{Klartag3},
and also on careful arguments similar to those in classical proofs of
the Berry-Esseen theorem.  It is not clear whether the Stein's method
approach can achieve these optimal error bounds.

\subsection*{Acknowledgements}
The author thanks S.\ Chatterjee and E.\ Meckes for showing him an
early version of \cite{CM} and E.\ Meckes for useful discussions.

\end{document}